\documentstyle[amsfonts,12pt,thmsa,sw20lart]{article}

\input{tcilatex}

\begin{document}

\title{Quantum It\^{o} B*-algebras, their\\
Classification and Decomposition.}
\author{V. P. Belavkin. \\
Mathematics Department, University of Nottingham,\\
NG7 2RD, UK.}
\date{July 20, 1997\\
Published in: {\it Quantum Probability} {\bf 43} 63--70 Banach Center
Publications, Warszawa 1998.}
\maketitle

\begin{abstract}
A simple axiomatic characterization of the general (infinite dimensional,
noncommutative) It\^{o} algebra is given and a pseudo-Euclidean fundamental
representation for such algebra is described. The notion of It\^{o}
B*-algebra, generalizing the C*-algebra is defined to include the Banach
infinite dimensional It\^{o} algebras of quantum Brownian and quantum L\'{e}%
vy motion, and the B*-algebras of vacuum and thermal quantum noise are
characterized. It is proved that every It\^{o} algebra is canonically
decomposed into the orthogonal sum of quantum Brownian (Wiener) algebra and
quantum L\'{e}vy (Poisson) algebra. In particular, every quantum thermal
noise is the orthogonal sum of a quantum Wiener noise and a quantum Poisson
noise as it is stated by the L\'{e}vy-Khinchin Theorem in the classical case
.
\end{abstract}

\section{Introduction: It\^{o} $\star $-algebras}

The classical stochastic calculus developed by It\^{o}, and its quantum
stochastic analog, given by Hudson and Parthasarathy (HP) in \cite{1}, can
be unified in the $\star $-algebraic approach to the operator integration in
Fock space \cite{2}, in which the classical and quantum calculi become
represented as two extreme commutative and completely noncommutative cases
of a generalized It\^{o} calculus.

On the axiomatic level the generalized It\^{o} algebra was defined in \cite%
{3} by a family of stochastic differentials ${\rm d}\Lambda \left(
t,a\right) $, $a\in {\frak a}$ with respect to $t\in {\Bbb R}_{+}$,
generating an associative, but in general noncommutative $\star $--algebra: $%
{\rm d}\Lambda \left( a\right) ^{\dagger }={\rm d}\Lambda \left( a^{\star
}\right) ,$ 
\begin{equation}
{\rm d}\Lambda \left( a\right) ^{\dagger }{\rm d}\Lambda \left( a\right) =%
{\rm d}\Lambda \left( a^{\star }a\right) ,\quad \sum \lambda _i{\rm d}%
\Lambda \left( a_i\right) ={\rm d}\Lambda \left( \sum \lambda _ia_i\right)
\label{0.1}
\end{equation}
with given mean values $\langle {\rm d}\Lambda \left( t,a\right) \rangle
=l\left( a\right) {\rm d}t$, absolutely continuous with respect to ${\rm d}t=%
{\rm d}\Lambda \left( t,\theta \right) $. Here the indexing $\star $%
-semigroup ${\frak a}$ is extended to a complex, in general infinite
dimensional space, parametrizing the It\^{o} algebra, with an involution $%
a\mapsto a^{\star }\in {\frak a}$, $a^{\star \star }=a$ and death $\theta
=\theta ^{\star }\in {\frak a}$, a self-adjoint annihilator ${\frak a}\theta
=0$, corresponding to ${\rm d}\Lambda \left( t,{\frak a}\right) {\rm d}t=0$,
and $l:{\frak a}\rightarrow {\Bbb C}$ is a positive $l\left( a^{\star
}a\right) \geq 0$ linear $*$-functional, $l\left( a^{\star }\right) =l\left(
a\right) ^{*}$, normalized as $l\left( \theta \right) =1$ due to the
determinism $\left\langle {\rm d}t\right\rangle ={\rm d}t$.

The functional $l$ defines the fundamental representation $%
\mbox{\boldmath$i$}:a\mapsto \mbox{\boldmath$a$}=\left( a_\nu ^\mu \right)
_{\nu =+,\bullet }^{\mu =-,\bullet }$ of ${\frak a}$ in terms of the
quadruples 
\begin{equation}
a_{\bullet }^{\bullet }=i\left( a\right) ,\quad a_{+}^{\bullet }=k\left(
a\right) ,\quad a_{\bullet }^{-}=k^{*}\left( a\right) ,\quad
a_{+}^{-}=l\left( a\right) ,  \label{0.2}
\end{equation}
where $i\left( a\right) =i\left( a^{\star }\right) ^{\dagger }$ is the GNS
representation $k\left( ab\right) =i\left( a\right) k\left( b\right) $ of $%
{\frak a}$ in the Hilbert space ${\cal H}\ni k\left( b\right) $, $b\in 
{\frak a}$ of the Kolmogorov decomposition $l\left( a^{\star }b\right)
=k\left( a\right) ^{\dagger }k\left( b\right) $, and $k^{*}\left( a\right)
=k\left( a^{\star }\right) ^{\dagger }$. In these tesor notations $a_\nu
^\mu =i_\nu ^\mu \left( a\right) $, $\mu =-,\bullet $, $\nu =+,\bullet $ the
composition properties $i\left( a\right) i\left( b\right) =i\left( ab\right) 
$, $k^{*}\left( a\right) i\left( b\right) =k^{*}\left( ab\right) $ together
with the already mentioned $i\left( a\right) k\left( b\right) =k\left(
ab\right) $, $k^{*}\left( a\right) k\left( b\right) =l\left( ab\right) ,$
can be written in the form of the multiplicativity $\mbox{\boldmath{i}}%
\left( ab\right) =\mbox{\boldmath$a$}\bullet \mbox{\boldmath$b$}$ with
respect to the convolution 
\begin{equation}
\mbox{\boldmath$a$}\bullet \mbox{\boldmath$b$}=\left( a_{\bullet }^\mu b_\nu
^{\bullet }\right) _{\nu =+,\bullet }^{\mu =-,\bullet }  \label{0.4}
\end{equation}
of the components $a_\nu $ and $b^\mu $ over the common index values $\mu
=\bullet =\nu $ only. One can also use the convenienve $a_{-}^\mu =0=a_\nu
^{+}$ of the tensor notations (\ref{0.2}) to extend the quadruples $%
\mbox{\boldmath$a$}=\mbox{\boldmath$i$}\left( a\right) $ to the matrices $%
{\bf a}=\left[ a_\nu ^\mu \right] _{\nu =-,\bullet ,+}^{\mu =-,\bullet ,+}$,
in which (\ref{0.4}) is simply given as ${\bf i}\left( ab\right) ={\bf ab}$
in terms of the usual product of the matrices ${\bf a}={\bf i}\left(
a\right) $ and ${\bf b}={\bf i}\left( b\right) $. But the involution $%
a\mapsto a^{\star }$, which is represented in terms of the quadruples as the
Hermitian conjugation $\mbox{\boldmath$i$}\left( a^{\star }\right) =\left(
a_{-\mu }^{-\nu \dagger }\right) _{\nu =+,\bullet }^{\mu =-,\bullet }$,
where $-(-)=+$, $-\bullet =\bullet $, $-(+)=-$, is given by the adjoint
matrix ${\bf ga}^{\dagger }{\bf g}={\bf i}\left( a^{\star }\right) $ w.r.t.
the pseudo-Hilbert (complex Minkowski) metrics ${\bf g}=\left[ \delta _{-\nu
}^\mu \right] _{\nu =-,\bullet ,+}^{\mu =-,\bullet ,+}$ as it was noted in 
\cite{2,3}.

It was proved in \cite{3,4} that any (classical or quantum) stochastic noise
described by a process $t\in {\Bbb R}_{+}\mapsto \Lambda \left( t,a\right)
,a\in {\frak a}$ with independent increments ${\rm d}\Lambda \left(
t,a\right) =\Lambda \left( t+{\rm d}t,a\right) -\Lambda \left( t,a\right) $,
forming an It\^{o} $*$-algebra, can be represented in the (symmetric) Fock
space ${\frak F}$ over the space of ${\cal H}$-valued square-integrable
functions on ${\Bbb R}_{+}$ as $\Lambda \left( t,a\right) =a_\nu ^\mu
\Lambda _\mu ^\nu \left( t\right) $. Here 
\begin{equation}
a_\nu ^\mu \Lambda _\mu ^\nu \left( t\right) =a_{\bullet }^{\bullet }\Lambda
_{\bullet }^{\bullet }\left( t\right) +a_{+}^{\bullet }\Lambda _{\bullet
}^{+}\left( t\right) +a_{\bullet }^{-}\Lambda _{-}^{\bullet }\left( t\right)
+a_{+}^{-}\Lambda _{-}^{+}\left( t\right) ,  \label{0.3}
\end{equation}
is the canonical decomposition of $\Lambda $ into the exchange $\Lambda
_{\bullet }^{\bullet }$, creation $\Lambda _{\bullet }^{+}$, annihilation $%
\Lambda _{-}^{\bullet }$ and preservation (time) $\Lambda _{-}^{+}=t{\rm I}$
operator-valued processes of the HP quantum stochastic calculus, having the
mean values $\left\langle \Lambda _\mu ^\nu \left( t\right) \right\rangle
=t\delta _{+}^\nu \delta _\mu ^{-}$ with respect to the vacuum state in $%
{\frak F}$. Thus the parametrizing algebra ${\frak a}$ of the processes with
independent increments can be always represented by a $\star $-subalgebra $%
\left\{ \mbox{\boldmath$i$}\left( a\right) |a\in {\frak a}\right\} $ of the
non-unital algebra ${\frak b}\left( {\cal H}\right) $ of all quadruples $%
\mbox{\boldmath$a$}=\left( a_\nu ^\mu \right) _{\nu =+,\bullet }^{\mu
=-,\bullet }$ with the product (\ref{0.4}) and the unique death, given by
the quadruple ${\bf \theta }=\left( \delta _{-}^\mu \delta _\nu ^{+}\right)
_{\nu =+,\bullet }^{\mu =-,\bullet }$ having the only nonzero element $%
\theta _{+}^{-}=1$. Here $a_\nu ^\mu :{\cal H}_\nu \rightarrow {\cal H}_\mu $
are the linear operators on ${\cal H}_{\bullet }={\cal H}$, ${\cal H}_{+}=%
{\Bbb C=}{\cal H}_{-}$, having the adjoints $a_\nu ^{\mu \dagger }:{\cal H}%
_\mu \rightarrow {\cal H}_\nu $, which define the natural involution $%
\mbox{\boldmath$a$}\mapsto \mbox{\boldmath$a$}^{\star }$ on ${\frak b}\left( 
{\cal H}\right) $, given by $a_{-\nu }^{\star \mu }=a_{-\mu }^{\nu \dagger }$%
.

The simplest, one-dimensional It\^{o} algebra generated by the nonstochastic
differential of the deterministic process $\Lambda \left( t,\theta \right) =t%
{\rm I}$ corresponds to the Newton differential calculus $\left( {\rm d}%
t\right) ^2=0$. It is described by the smallest It\^{o} algebra ${\frak a}=%
{\Bbb C}\theta $ with $l\left( a\right) =\alpha \in {\Bbb C}$ for $a=\alpha
\theta $ and the nilpotent multiplication $\alpha ^{\star }\alpha =0$ with $%
\alpha ^{\star }=\bar{\alpha}$ representing the classical non-stochastic
calculus in ${\cal H}=\left\{ 0\right\} $.

The classical stochastic calculus ${\rm d}w_i{\rm d}w_k=\delta _k^i{\rm d}t$
for the It\^{o} differentials ${\rm d}w_i$, $i=1,\ldots ,d$ of the standard
Wiener processes $w_i\left( t\right) $ generates the $d+1$ dimensional It%
\^{o} algebra as the second order nilpotent algebra ${\frak a}$ of pairs $%
a=\left( \alpha ,\eta \right) $, $\eta \in {\Bbb C}^d$ with $\theta =\left(
1,0\right) $ and $a^{\star }a=\left( \left\| \eta \right\| ^2,0\right) $ for 
$a^{\star }=\left( \bar{\alpha},\eta ^{\star }\right) $, where $\eta ^{\star
}=\left( \bar{\eta}^i\right) $ is the complex conjugation of $\eta =\left(
\eta ^i\right) $ and $\left\| \eta \right\| ^2=\sum_{i=1}^d\left| \eta
^i\right| ^2$. The GNS representation associated with $l\left( a\right)
=\alpha $ is given by $a_{+}^{-}=\alpha ,\quad a_i^{-}=\eta ^i=a_{+}^i,\quad
a_{\bullet }^{\bullet }=0$ in ${\cal H}={\Bbb C}^d$, and the operator
representation of ${\frak a}$ in Fock space is defined by the forward
differentials of $\Lambda \left( t,a\right) =\alpha t{\rm I}%
+\sum_{i=1}^d\eta ^i{\rm W}_i\left( t\right) $, where ${\rm W}_i\left(
t\right) =\Lambda _{-}^i\left( t\right) +\Lambda _i^{+}\left( t\right) $ are
the operator representations of the processes $w_i\left( t\right) $ with
respect to the vacuum state in ${\frak F}$.

The unital $\star $-algebra of complex vectors $\zeta =\left( \zeta
^i\right) $ with component-vise multiplication $\zeta ^{\star }\zeta =\left(
\left| \zeta ^i\right| ^2\right) =\left| \zeta \right| ^2$ can be embedded
into the $d+1$-dimensional It\^{o} algebra ${\frak a}$ of the pairs $%
a=\left( \alpha ,\zeta \right) $, $\zeta \in {\Bbb C}^d$ with $\theta
=\left( 1,0\right) $ and $a^{\star }a=\left( \left\| \zeta \right\|
^2,\left| \zeta \right| ^2\right) $ by $\alpha =\sum_{i=1}^d\zeta ^i$. The
GNS representation of ${\frak a}$, associated with $l\left( a\right) =\alpha 
$, is given by $a_{+}^{-}=\alpha $, $a_i^{-}=\zeta ^i=a_{+}^i$, $%
a_k^i=\delta _k^i\zeta ^i$, and the Fock space representation is defined by $%
\Lambda \left( t,a\right) =\alpha t{\rm I}+\sum_{i=1}^d\zeta ^i{\rm M}%
_i\left( t\right) $, where ${\rm M}_i\left( t\right) =\Lambda _i^i\left(
t\right) +{\rm W}_i\left( t\right) $ are the operator representations in $%
{\frak F}$ of the standard Poisson processes $n_i\left( t\right) =m_i\left(
t\right) +t$, compensated by their vacuum mean values $\left\langle
n_i\left( t\right) \right\rangle =t$, with ${\rm d}n_i{\rm d}n_k=\delta _k^i%
{\rm d}n_k$.

Note that any two-dimensional It\^{o} $\star $-algebra ${\frak a}$ is
commutative as $\theta a=0=a\theta $ for any other element $a\neq \theta $
of a basis $\left\{ a,\theta \right\} $ in ${\frak a}$. Moreover, each such
algebra is either of the Wiener or of the Poisson type as it is either
second order nilpotent, or contains a unital one-dimensional subalgebra.
Thus our results on the classification of It\^{o} $\star $-algebras will be
nontrivial only in the higher dimensions of ${\frak a}$. The well known
Levy-Khinchin classification of the classical noise can be reformulated in
purely algebraic terms as the decomposability of any commutative It\^{o}
algebra into the Newton, Wiener (Brownian) and Poisson (L\'{e}vy) orthogonal
components. In the general case we shall show that every It\^{o} $\star $%
-algebra is also decomposable into the Newton, a quantum Brownian, and a
quantum L\'{e}vy orthogonal components.

Let us first consider the following class of strictly non-commutative vacuum
It\^{o} algebras which have no classical analogs, and then a more ordinary
class of thermal It\^{o} algebras as they were defined in \cite{3,4} by a
generalization of classical Wiener and Poisson noise algebras.

\section{It\^{o} B*-algebra of vacuum noise}

Let ${\frak h}$ be a Hilbert space of ket-vectors $x$ with scalar product $%
\left\langle x|x\right\rangle $ and ${\cal A}\subseteq {\cal B}\left( {\frak %
h}\right) $ be a C*-algebra, represented on ${\frak h}$ by the operators $%
{\cal A}\ni A:x\mapsto Ax$ with $\left\langle A^{\dagger }x|x\right\rangle
=\left\langle x|Ax\right\rangle $. We denote by ${\frak h}^{*}$ the dual
Hilbert space of bra-vectors $y=z^{*}$, $z\in {\frak h}$ with the scalar
product $\left\langle y|x^{*}\right\rangle =yx=\left\langle
y^{*}|x\right\rangle $ given by inverting anti-linear isomorphism ${\frak h}%
^{*}\ni y\mapsto y^{*}\in {\frak h}$, and the dual representation of ${\cal A%
}$ as the right representation $A^{\prime }:y\mapsto yA$, $y\in {\frak h}%
^{*} $, given by $\left( yA\right) x=y\left( Ax\right) $ such that $%
\left\langle yA^{\dagger }|y\right\rangle =\left\langle y|yA\right\rangle $
on ${\frak h}^{*}$. Then the direct sum ${\cal K}={\frak h}\oplus {\frak h}%
^{*}$ becomes a two-sided ${\cal A}$-module 
\begin{equation}
A\xi =Ax,\quad \xi A=yA,\quad \forall \xi =x\oplus y,  \label{1.1}
\end{equation}
with the flip-involution $\left( x\oplus y\right) ^{\star }=y^{*}\oplus
x^{*} $ and two scalar products 
\begin{equation}
\left\langle x\oplus y^{\prime }|x^{\prime }\oplus y\right\rangle
_{+}=\left\langle x|x^{\prime }\right\rangle ,\quad \left\langle x\oplus
y^{\prime }|x^{\prime }\oplus y\right\rangle ^{-}=\left\langle y^{\prime
}|y\right\rangle .  \label{1.2}
\end{equation}

The space ${\frak a}={\Bbb C}\oplus {\cal K}\oplus {\cal A}$ of triples $%
a=\left( \alpha ,\xi ,A\right) $ becomes an It\^{o} $\star $-algebra with
respect to the non-commutative product 
\begin{equation}
a^{\star }a=\left( \left\langle \xi |\xi \right\rangle _{+},\xi ^{\star
}A+A^{\dagger }\xi ,A^{\dagger }A\right) ,\quad aa^{\star }=\left(
\left\langle \xi |\xi \right\rangle ^{-},A\xi ^{\star }+\xi A^{\dagger
},AA^{\dagger }\right) ,  \label{1.3}
\end{equation}
where $\left( \alpha ,\xi ,A\right) ^{\star }=\left( \bar{\alpha},\xi
^{\star },A^{\dagger }\right) $, with death $\theta =\left( 1,0,0\right) $
and $l\left( \alpha ,\xi ,A\right) =\alpha $. Obviously $a^{\star }a\neq
aa^{\star }$ if $\left\| \xi \right\| _{+}=\left\| x\right\| \neq $ $\left\|
y\right\| =\left\| \xi \right\| ^{-}$ even if the operator algebra ${\cal A}$
is commutative, $A^{\dagger }A=AA^{\dagger }$. It is separated by four
semi-norms 
\begin{equation}
\left\| a\right\| =\left\| A\right\| ,\left\| a\right\| _{+}=\left\|
x\right\| ,\left\| a\right\| ^{-}=\left\| y\right\| ,\left\| a\right\|
_{+}^{-}=\left| \alpha \right| ,  \label{1.4}
\end{equation}
and is jointly complete as $a=\left( \alpha ,x\oplus y,A\right) \in {\frak a}
$ have independent components from the Banach spaces ${\Bbb C}$, ${\frak h}$%
, ${\frak h}^{*}$and ${\cal A}$.

We shall call such Banach It\^{o} algebra the vacuum algebra as $l\left(
a^{\star }a\right) =0$ for any $a\in {\frak a}$ with $\xi \in {\frak h}^{*}$
(the Hudson-Parthasarathy algebra ${\frak a}={\frak b}\left( {\frak h}%
\right) $ if ${\cal A}={\cal B}\left( {\frak h}\right) $). Every closed It%
\^{o} subalgebra ${\frak a}\subseteq {\frak b}\left( {\cal H}\right) $ of
the HP algebra ${\frak b}\left( {\cal H}\right) $ equipped with four norms (%
\ref{1.4}) on a Hilbert space ${\cal H}$ is called the operator It\^{o}
B*-algebra.

If the algebra ${\cal A}$ is completely degenerated on ${\frak h}$, ${\cal A}%
=\left\{ 0\right\} $, the It\^{o} algebra ${\frak a}$ is nilpotent of second
order, and contains only the two-dimensional subalgebras of Wiener type $%
{\frak b}={\Bbb C}\oplus {\Bbb C}\oplus \left\{ 0\right\} $ generated by an $%
a=\left( \alpha ,x\oplus y,0\right) $ with $\left\| x\right\| =\left\|
y\right\| $. Every closed It\^{o} subalgebra ${\frak b}\subseteq {\frak a}$
of the HP B*-algebra ${\frak a}={\frak b}\left( {\cal H}\right) $ is called
the B*-It\^{o} algebra of a vacuum Brownian motion if it is defined by a $%
\star $-invariant direct sum ${\cal G}={\cal G}_{+}\oplus {\cal G}%
^{-}\subseteq {\cal K}$ given by a Hilbert subspace ${\cal G}_{+}\subseteq 
{\cal H}$, ${\cal G}^{-}={\cal G}_{+}^{*}$ and ${\cal A}=\left\{ 0\right\} $.

In the case $I\in {\cal A}$ the algebra ${\cal A}$ is not degenerated and
contains also the vacuum Poisson subalgebra ${\Bbb C}\oplus \left\{
0\right\} \oplus {\Bbb C}I$ of the total quantum number on ${\frak h}$, and
other Poisson two-dimensional subalgebras, generated by $a=\left( \alpha
,x\oplus y,I\right) $ with $y=e^{i\phi }x^{*}$. We shall call a closed It%
\^{o} subalgebra ${\frak c}\subseteq {\frak a}$ of the HP B*-algebra ${\frak %
a}={\frak b}\left( {\cal H}\right) $ the B*-algebra of a vacuum L\'{e}vy
motion if it is given by a direct sum ${\cal E}={\cal E}_{+}\oplus {\cal E}%
^{-}\subseteq {\cal K}$ with ${\cal E}^{\_}={\cal E}_{+}^{*}$ and a $*$%
-subalgebra ${\cal A}\subseteq {\cal B}\left( {\cal H}\right) $
nondegenerated on the subspace ${\cal E}_{+}\subseteq {\cal H}$.

We shall see that the general vacuum B*-algebra ${\frak a}$, which is
characterized by the condition ${\frak n}_{+}^{\perp }={\frak n}^{-}$, where 
\begin{equation}
{\frak n}^{-}=\left\{ b\in {\frak a}:\left\| b\right\| ^{-}=0\right\} ,\quad 
{\frak n}_{+}=\left\{ c\in {\frak a}:\left\| c\right\| _{+}=0\right\} ,
\label{1.5}
\end{equation}
and ${\frak n}_{+}^{\perp }$ is the right orthogonal complement to ${\frak n}%
_{+}$, can be represented as an operator vacuum It\^{o} B*-algebra.

\begin{theorem}
Every vacuum B*-algebra can be decomposed into an orthogonal sum ${\frak a}=%
{\frak b}+{\frak c}$, ${\frak bc}=\left\{ 0\right\} $ of the Brownian vacuum
B*-algebra ${\frak b}$ and the L\'{e}vy vacuum B*-algebra ${\frak c}$.
\end{theorem}

Proof. This decomposition is uniquely defined for all $a=\left( \alpha ,\xi
,A\right) $ by $a=\alpha \theta +b+c$, with $b=\left( 0,\eta ,0\right) $, $%
c=\left( 0,\zeta ,A\right) $, $\eta =Px\oplus yP\in {\cal G}$, $\zeta =\xi
-\eta \in {\cal E}$, where $P=P^{\dagger }$ is the maximal projector in $%
{\frak h}$, for which ${\cal A}P=\left\{ 0\right\} $, ${\cal G}_{+}=P{\frak h%
}$, and ${\cal E}_{+}{\cal =G}_{+}^{\perp }$. End of proof.

\section{It\^{o} B*-algebra of thermal noise}

Let ${\cal D}$ be a left Tomita $\star $-algebra \cite{5} with respect to a
Hilbert norm $\left\| \xi \right\| _{+}=0\Rightarrow \xi =0$, and thus a
right pre-Hilbert $\star $-algebra with respect to $\left\| \xi \right\|
^{-}=\left\| \xi ^{\star }\right\| _{+}$. This means that ${\cal D}$ is a
complex pre-Hilbert space with continuous left (right) multiplications $%
C:\zeta \mapsto \xi \zeta $ ($C^{\prime }:\eta \mapsto \eta \xi $) w.r.t. $%
\left\| \cdot \right\| _{+}$ (w.r.t. $\left\| \cdot \right\| ^{-}$) of the
elements $\zeta ,\eta \in {\cal D}$ respectively, defined by an associative
product in ${\cal D}$, and the involution ${\cal D}\ni \xi \mapsto \xi
^{\star }\in {\cal D}$ such that 
\begin{equation}
\left\langle \eta \zeta ^{\star }|\xi \right\rangle ^{-}=\left\langle \eta
|\xi \zeta \right\rangle ^{-},\quad \left\langle \eta ^{\star }\zeta |\xi
\right\rangle _{+}=\left\langle \zeta |\eta \xi \right\rangle _{+}\quad
\forall \xi ,\zeta ,\eta \in {\cal D},  \label{2.1}
\end{equation}
\begin{equation}
\left\langle \eta |\xi ^{\star }\right\rangle ^{-}=\left\langle \xi |\eta
^{\sharp }\right\rangle ^{-},\quad \left\langle \zeta |\xi ^{\star
}\right\rangle _{+}=\left\langle \xi |\zeta ^{\flat }\right\rangle _{+}\quad
\forall \eta \in {\cal D}^{-},\zeta \in {\cal D}_{+}.  \label{2.2}
\end{equation}
Here $\left\langle \eta |\zeta ^{\star }\right\rangle ^{-}=\left\langle \eta
^{\star }|\zeta \right\rangle _{+}$ is the right scalar product, ${\cal D}%
_{+}={\cal D}_{+}^{\flat }$ is a dense domain for the left adjoint
involution $\zeta \mapsto \zeta ^{\flat }$, $\zeta ^{\flat \flat }=\zeta $,
and ${\cal D}^{-}={\cal D}_{+}^{\star }$ is the invariant domain for the
right adjoint involution $\eta \mapsto \eta ^{\sharp }$, $\left( \eta
^{\sharp }\eta \right) ^{\sharp }=\eta ^{\sharp }\eta $ such that $\zeta
^{\flat \star }=\zeta ^{\star \natural }$, $\eta ^{\sharp \star }=\eta
^{\star \flat }$.

Since the adjoint operators $C^{\dagger }\zeta =\xi ^{\star }\zeta $, $\eta
C^{\dagger }=\eta \xi ^{\star }$ are also given by the multiplications, they
are bounded: 
\begin{equation}
\left\| \xi \right\| =\sup \left\{ \left\| \xi \zeta \right\| _{+}:\left\|
\zeta \right\| _{+}\leq 1\right\} =\sup \left\{ \left\| \eta \xi \right\|
^{-}:\left\| \eta \right\| ^{-}\leq 1\right\} <\infty .  \label{2.3}
\end{equation}
Note that we do not require the sub-space ${\cal DD}\subseteq {\cal D}$ of
all products $\eta \zeta $, $\eta ,\zeta \in {\cal D}$ to be dense in ${\cal %
D}$ w.r.t. any of two Hilbert norms on ${\cal D}$, but it is always dense
w.r.t. the operator semi-norm (\ref{2.3}) on ${\cal D}$. Hence the operator $%
\dagger $-algebra ${\cal C}=\left\{ C:{\cal D}\ni \zeta \mapsto \xi \zeta
|\xi \in {\cal D}\right\} $ w.r.t. the left scalar product, which is also
represented on the ${\cal D}\ni \eta $ equipped with $\left\langle \cdot
|\cdot \right\rangle ^{-}$ by the right multiplications $\eta C=\eta \xi $, $%
\xi \in {\cal D}$, can be degenerated on ${\cal D}$.

Thus the direct sum ${\frak a}={\Bbb C}\oplus {\cal D}$ of pairs $a=\left(
\alpha ,\xi \right) $ becomes an It\^{o} $\star $-algebra with the product 
\begin{equation}
a^{\star }a=\left( \left\langle \xi |\xi \right\rangle _{+},\xi ^{\star }\xi
\right) ,\quad aa^{\star }=\left( \left\langle \xi |\xi \right\rangle
^{-},\xi \xi ^{\star }\right) ,  \label{2.4}
\end{equation}
where $\left( \alpha ,\xi \right) ^{\star }=\left( \bar{\alpha},\xi ^{\star
}\right) $, with death $\theta =\left( 1,0\right) $ and $l\left( \alpha ,\xi
\right) =\alpha $. Obviously $a^{\star }a\neq aa^{\star }$ if the involution 
$a\mapsto a^{\star }$ is not isometric w.r.t. any of two Hilbert norms even
if the algebra ${\cal D}$ is commutative. It is a Banach algebra if its
normed $\star $-algebra ${\cal D}$ is complete jointly w.r.t. to three
norms, and is separated by the semi-norms 
\begin{equation}
\left\| a\right\| =\left\| \xi \right\| ,\left\| a\right\| ^{-}=\left\| \xi
\right\| ^{-},\left\| a\right\| _{+}=\left\| \xi \right\| _{+},\left\|
a\right\| _{+}^{-}=\left| \zeta \right| .  \label{2.5}
\end{equation}

We shall call such complete It\^{o} algebra the thermal B*-algebra as $%
l\left( a^{\star }a\right) =\left\| \xi \right\| _{+}^2\neq 0$ for any $a\in 
{\frak a}$ with $\xi \neq 0$. If $\zeta \eta =0$ for any $\zeta ,\eta \in 
{\cal D}$, it is the It\^{o} B*-algebra of thermal Brownian motion. A
thermal B*-subalgebra ${\frak b}\subseteq {\frak a}$ with such trivial
product is given by any involutive pre-Hilbert $\star $-invariant two-normed
subspace ${\cal G}\subseteq {\cal D}$ which is closed w.r.t. the Hilbert sum 
$\left\langle \eta |\zeta \right\rangle ^{-}+\left\langle \xi |\zeta
\right\rangle _{+}$. We shall call such Brownian algebra ${\frak b}={\Bbb %
C\oplus }{\cal G}$ the quantum (if $\left\| \cdot \right\| _{+}\neq $ $%
\left\| \cdot \right\| ^{-}$) Wiener B*-algebra associated with the space $%
{\cal G}$.

In the opposite case, if ${\cal DD}=\left\{ \zeta \eta :\zeta ,\eta \in 
{\cal D}\right\} $ is dense in ${\cal D}$, it has nondegenerated operator
representation ${\cal C}$ on ${\cal D}$. Any closed involutive sub-algebra $%
{\cal E}\subseteq {\cal D}$ which is non-degenerated on ${\cal E}$ defines
an It\^{o} B*-algebra ${\frak c}={\Bbb C}\oplus {\cal E}$ of thermal L\'{e}%
vy motion. We shall call such It\^{o} algebra the quantum (if ${\cal E} $ is
non-commutative) Poisson B*-algebra.

In next section the general thermal algebra ${\frak a}$, characterized as an
It\^{o} B*-algebra with the condition of closability of the involution $%
a\mapsto a^{\star }$ w.r.t. any of left or right Hilbert semi-norms on $%
{\frak a}$, is represented as a concrete thermal B*-algebra ${\Bbb C}\oplus 
{\cal D}$.

\begin{theorem}
Every thermal It\^{o} B*-algebra is an orthogonal sum ${\frak a}={\frak b}+%
{\frak c}$, ${\frak bc}=\left\{ 0\right\} $ of the Wiener B*-algebra ${\frak %
b}$ and the Poisson B*-algebra ${\frak c}$.
\end{theorem}

Proof. The orthogonal decomposition $a=\alpha \theta +b+c$ for all $a=\left(
\alpha ,\xi \right) \in {\frak a}$, uniquely given by the decomposition $\xi
=\eta +\zeta $ w.r.t. any of two scalar products in ${\cal D}$, where $\eta
=P\xi =\xi P$ is the orthogonal projection onto ${\cal G}\perp {\cal DD}$
w.r.t. any of two Hilbert norms, and $\zeta =\xi -\eta $ .

Indeed, if $\xi \in {\cal D}$ is left orthogonal to ${\cal DD}$, then it is
also right orthogonal to ${\cal DD}$ and vice versa: 
\begin{eqnarray*}
\left\langle \eta \zeta ^{\star }|\xi \right\rangle ^{-} &=&\left\langle
\zeta \eta ^{\star }|\xi ^{\star }\right\rangle _{+}=\left\langle \xi |\eta
^{\sharp \star }\zeta ^{\flat }\right\rangle _{+}=0,\quad \forall \eta \in 
{\cal D}^{-},\zeta \in {\cal D}_{+}, \\
\left\langle \eta ^{\star }\zeta |\xi \right\rangle _{+} &=&\left\langle
\zeta ^{\star }\eta |\xi ^{\star }\right\rangle ^{-}=\left\langle \xi |\eta
^{\sharp }\zeta ^{\flat \star }\right\rangle ^{-}=0,\quad \forall \eta \in 
{\cal D}^{-},\zeta \in {\cal D}_{+}.
\end{eqnarray*}
From these and (\ref{2.1}) equations it follows that $\eta \xi =0=\xi \zeta $
for all $\zeta ,\eta \in {\cal D}$ if $\xi $ is (right or left) orthogonal
to ${\cal DD}$, and so $\left\| \xi \right\| =0$ for such $\xi $ and vice
versa. Thus the orthogonal subspace is the space ${\cal G}=\left\{ \xi \in 
{\cal D}:\left\| \xi \right\| =0\right\} $ is jointly complete w.r.t. only
two norms $\left\| \cdot \right\| ^{-},\left\| \cdot \right\| _{+}$ and can
be considered as the Hilbert subspace of pre-Hilbert space ${\cal D}$ with
isometric involution $\xi \mapsto \xi ^{\star }$ w.r.t. $\left\langle \eta
|\zeta \right\rangle ^{-}+\left\langle \zeta |\eta \right\rangle _{+}$.
Denoting by $P$ the orthogonal projector in ${\cal D} $ onto ${\cal G}%
\subseteq {\cal D}$ we obtain $P\xi =\eta =\xi P$ and $\zeta =\xi -\eta \in 
{\cal D}$ is in the closure ${\cal E}\subseteq {\cal D}$ of ${\cal DD}$. End
of proof.

\section{Decomposition of It\^{o} B*-algebras}

Now we shall consider the general case. Let ${\frak a}$ be an associative
infinite-dimensional complex algebra with involution $b^{\star }=a\in {\frak %
a},\forall b=a^{\star }$ which is defined by the properties

\[
\left( a^{\star }a\right) ^{\star }=a^{\star }a,\quad \left( \sum \lambda
_ib_i\right) ^{\star }=\sum \bar{\lambda}_ib_i^{\star },\quad \forall b_i\in 
{\frak a},\lambda _i\in {\Bbb C}. 
\]
We shall suppose that this algebra is a normed space with respect to four
semi-norms $\left\| \cdot \right\| _\nu ^\mu $, indexed as 
\begin{equation}
\left\| \cdot \right\| _{\bullet }^{\bullet }\equiv \left\| \cdot \right\|
,\quad \left\| \cdot \right\| _{+}^{\bullet }\equiv \left\| \cdot \right\|
_{+},\quad \left\| \cdot \right\| _{\bullet }^{-}\equiv \left\| \cdot
\right\| ^{-},\quad \left\| \cdot \right\| _{+}^{-}  \label{3.1}
\end{equation}
by $\mu =-,\bullet $, $\nu =+,\bullet $, satisfying the following conditions

\begin{equation}
\left\| b^{\star }\right\| =\left\| b\right\| ,\quad \left\| b^{\star
}\right\| _{+}=\left\| b\right\| ^{-},\quad \left\| b^{\star }\right\|
_{+}^{-}=\left\| b\right\| _{+}^{-},  \label{3.2}
\end{equation}
\begin{equation}
\left( \left\| ac\right\| _\nu ^\mu \leq \left\| a\right\| _{\bullet }^\mu
\left\| c\right\| _\nu ^{\bullet },\right) _{\nu =+,\bullet }^{\mu
=-,\bullet }\quad \forall \,a,b,c\in {\frak a}.  \label{3.3}
\end{equation}
Thus the semi-norms (\ref{3.1}) separate ${\frak a}$ in the sense 
\[
\left\| a\right\| =\left\| a\right\| _{+}=\left\| a\right\| ^{-}=\left\|
a\right\| _{+}^{-}=0\quad \Rightarrow \quad a=0, 
\]
and the product $\left( a,c\right) \mapsto ac$ with involution $\star $ are
uniformly continuous in the induced topology due to (\ref{3.3}).

If ${\frak a}$ is a $*$-algebra equipped with a linear positive $*$%
-functional $l$ such that 
\[
l\left( a\right) =l\left( ac\right) =l\left( c^{\prime }a\right) =l\left(
c^{\prime }ac\right) =0,\forall c,c^{\prime }\in {\frak a\quad }\Rightarrow
\quad a=0, 
\]
and it is bounded with respect to $l$ in the sense 
\begin{equation}
\left\| a\right\| =\sup \left\{ \left\| c^{\prime }ac\right\|
_{+}^{-}/\left\| c^{\prime }\right\| ^{-}\left\| c\right\| _{+}:c^{\prime
},c\in {\frak a}\right\} <\infty \quad \forall a\in {\frak a,}  \label{3.4}
\end{equation}
where $\left\| a\right\| _{+}^{-}=\left| l\left( a\right) \right| ,\quad
\left\| a\right\| ^{-}=\left( l\left( aa^{\star }\right) \right)
^{1/2},\quad \left\| a\right\| _{+}=\left( l\left( a^{\star }a\right)
\right) ^{1/2}$, then it is four-normed in the above sense. The defined by $%
l $ semi-norms $\left\| \cdot \right\| _\nu ^\mu $ are obviously separating,
satisfy the inequalities (\ref{3.3}), and they satisfy also the $\star $%
-equalities of the following definition

\begin{definition}
An associative four-normed $\star $-algebra ${\frak a}$ is called B*-algebra
if it is complete in the uniform topology, induced by the semi-norms $\left(
\left\| a\right\| _\nu ^\mu \right) _{\nu =+,\bullet }^{\mu =-,\bullet }$,
satisfying the following equalities 
\begin{equation}
\left\| a^{\star }a\right\| =\left\| a^{\star }\right\| \left\| a\right\|
,\quad \left\| a^{\star }a\right\| _{+}^{-}=\left\| a^{\star }\right\|
^{-}\left\| a\right\| _{+}\quad \forall a\in {\frak a}..  \label{3.5}
\end{equation}
The It\^{o} B*-algebra is a B*-algebra with self-adjoint annihilator $\theta
=\theta ^{\star }$, $a\theta =0=\theta a$, $\forall a\in {\frak a}$ called
death for ${\frak a}$, and the semi-norms (\ref{3.1}) given by a linear
positive $*$-functional $l\left( a^{\star }\right) =l\left( a\right) ^{*}$, $%
l\left( a^{\star }a\right) \geq 0,\forall a\in {\frak a}$ normalized as $%
l\left( \theta \right) =1$.
\end{definition}

Obviously, any C*-algebra can be considered as a B*-algebra in the above
sense with three trivial semi-norms $\left\| a\right\| _{+}^{-}=\left\|
a\right\| ^{-}=\left\| a\right\| _{+}=0,\forall a\in {\frak a}$. Moreover,
as it follows from the inequalities (\ref{3.3}) for $c=1$, every unital
B*-algebra is a C*-algebra, the three nontrivial semi-norms on which might
be given by a state $l$, normalized as $l\left( 1\right) =\left\| 1\right\|
_{+}^{-}$. However, if a B*-algebra ${\frak a}$ contains only approximative
identity $e_i\nearrow 1$, and $\left\| e_i\right\| _{+}^{-}\longrightarrow
\infty $, it is a proper dense sub-algebra of its C*-algebraic completion
w.r.t. the norm $\left\| \cdot \right\| $.

Note that B*-algebra is an old name for the C*-algebra, so the use of the
obsolete term in a slightly more general sense is not contradictive.
Moreover, the term It\^{o} B*-algebra, or more appropriately, B*-It\^{o}
algebra will never make a confusion, as there is no It\^{o} algebra which is
simultaneously a C*-algebra. Indeed, suppose that an It\^{o} B*-algebra is a
C*-algebra. Then the functional $l$ is normalized to $1=l\left( \theta
\right) $ on $\theta =0$, the only annihilator $\theta \in {\frak a}$ in a
C*-algebra, as from $\left\| \theta \right\| ^2=\left\| \theta \theta
\right\| =\left\| 0\right\| =0$ it follows that $\theta =0$ in any
C*-algebra. But the condition $l\left( 0\right) =1$ contradicts to the
linearity of $l$.

It was proved in \cite{3,4} as mentioned already in the introduction, every
It\^{o} B*-algebra ${\frak a}$ is algebraically and isometrically isomorphic
to a closed subalgebra of the vacuum\ It\^{o} B*-algebra ${\frak b}\left( 
{\cal H}\right) $ of a Hilbert space ${\cal H}$ by $\mbox{\boldmath$i$}%
\left( a\right) =\left( \alpha ,\xi ,A\right) $, where 
\[
\alpha =l\left( a\right) ,\quad \xi =k\left( a\right) \oplus k^{*}\left(
a\right) ,\quad A=i\left( a\right) .\quad 
\]
It is a vacuum B*-algebra: ${\cal K}={\cal K}_{+}\oplus {\cal K}^{-}$, where 
${\cal K}_{+}=k\left( {\frak a}\right) $, ${\cal K}^{-}=k^{*}\left( {\frak a}%
\right) $ iff from orthogonality of $k^{*}\left( a\right) $ to all $%
k^{*}\left( c\right) $ with $k\left( c\right) =0$ it follows $k^{*}\left(
a\right) =0$. This can be expressed as ${\frak k}_{+}={\frak n}^{-}$ or $%
{\frak k}^{-}={\frak n}_{+}$ in terms of the right orthogonal complement $%
{\frak k}_{+}={\frak n}_{+}^{\perp }$ and 
\begin{equation}
{\frak k}^{-}=\left\{ a\in {\frak a}:\left\langle a|b\right\rangle
_{+}=0,\forall b\in {\frak n}^{-}\right\} ={\frak k}_{+}^{\star }
\label{3.6}
\end{equation}
and the ideals (\ref{1.5}), as ${\frak n}_{+}=k^{-1}\left( 0\right) $ is the
(left) null ideal for the map $k$ and ${\frak n}^{-}={\frak n}_{+}^{\star }$%
. It follows from the canonical construction of ${\cal K}^{-}={\cal K}%
_{+}^{*}$ as the quotient space ${\frak a}/{\frak n}^{-}$.

Due to the orthogonality of ${\frak k}_{+}$ and ${\frak k}^{-}$ in vacuum It%
\^{o} algebras, the involution $\star $ is never defined in ${\frak k}_{+}$
or in ${\frak k}^{-}$ except on the jointly null ideal ${\frak n}^{-}\cap 
{\frak n}_{+}$. The thermal B*-algebras have the trivial ideals ${\frak n}%
_{+}={\Bbb C}\theta ={\frak n}^{-}$, and so the involution $\star $ is
defined into ${\frak k}_{+}$ on the whole ${\frak k}_{+}={\frak a}={\frak k}%
^{-}$, and thus on the pre-Hilbert space ${\cal D}={\frak a}/{\Bbb C}\theta $
identified with $\left\{ \xi =a-l\left( a\right) \theta :a\in {\frak a}%
\right\} $, by $\xi ^{\star }=a^{\star }-l\left( a^{\star }\right) \theta $
. If this involution is left or right isometric (the case of tracial It\^{o}
algebras), it coincides with the left and right adjoint involutions (\ref%
{2.2}). If it is not isometric, it has densely defined left and right
adjoints iff it is left or right closable on ${\cal D}$. Thus the general
thermal B*-algebra can be characterized as an It\^{o} B*-algebra with
closable involution w.r.t. the Hilbert norms on ${\cal D}$.

In the general case the involution $\star $ has the range in ${\frak k}_{+}$
only if it is restricted to the invariant domain ${\frak d}={\frak k}%
^{-}\cap {\frak k}_{+}$. Thus we can say that it is closable in ${\frak k}%
_{+}$ if $\left\| a_n\right\| _{+}\longrightarrow 0\Rightarrow \left\|
a_n\right\| ^{-}\longrightarrow 0$ for any fundamental w.r.t. both Hilbert
semi-norms sequence $a_n\in {\frak d}$ and so it is also closable in ${\frak %
k}^{-}$. We shall call an It\^{o} B*-algebra ${\frak a}$ with a closable in
this sense involution the B*-Ito algebra, or the general Brownian algebra if
the semi-norm $\left\| \cdot \right\| =0$ is trivial on ${\frak a}$, and the
general L\'{e}vy algebra in the opposite case, when ${\frak aa}$ is dense in 
${\frak a}$ w.r.t. any of the Hilbert semi-norms $\left\| \cdot \right\|
^{-},\left\| \cdot \right\| _{+\,}$ (it is always dense w.r.t. the operator
semi-norm $\left\| \cdot \right\| $).

\begin{theorem}
Let ${\frak a}$ be an It\^{o} B*-algebra in which the restricted involution $%
\star :{\frak d}\rightarrow {\frak a}$ on the orthogonal complement ${\frak d%
}=\left( {\frak n}^{-}\cup {\frak n}_{+}\right) ^{\bot }$ to (\ref{1.5})
w.r.t. both Hilbert semi-norms is left (or right) closable. Then it is an
orthogonal sum ${\frak b}+{\frak c}$, ${\frak bc}=0$ of a quantum Brownian
B*-algebra ${\frak b}$ and a quantum L\'{e}vy B*-algebra ${\frak c}$.
\end{theorem}

Proof. Here we give the proof under a stronger than closability condition
when the quotient algebra ${\frak a}/{\frak n}$ with respect to the null $%
\star $-ideal ${\frak n}=\left\{ a\in {\frak a}:\left\| a\right\| =0\right\} 
$ has an identity, which defines the supporting ortho-projector $E=i\left(
e\right) $ for operator representation ${\cal A}=i\left( {\frak a}\right)
\simeq {\frak a}/{\frak n}$ on a Hilbert space ${\cal H}$. This means that
there exists an element $e=e^{\star }\in {\frak a}$ such that $i\left(
ae\right) =i\left( a\right) =i\left( ea\right) $ for all $a\in {\frak a}$, 
\[
k\left( aec\right) =i\left( ae\right) k\left( c\right) =k\left( ac\right)
,\quad k^{*}\left( aec\right) =k^{*}\left( c\right) i\left( ea\right)
=k^{*}\left( ac\right) , 
\]
and so $aec=ac-l\left( ac-aec\right) \theta $ for all $a,c\in {\frak a}.$ as
it is in the canonical representation (\ref{0.2}). We assume that $e$ is an
idempotent, otherwise it should be replaced by $e^2$.

Under this condition, which is automatically fulfilled in the finite
dimensional case, we can easily define the required orthogonal decomposition 
${\frak a}={\frak b}+{\frak c}$ by 
\[
a=b+c,\quad c=ae+ea-eae. 
\]
Here $b$ is an element of the quantum Brownian algebra ${\frak b}=\left\{
b\in {\frak a}:be=0\right\} \subseteq {\frak n}$ which is closed in ${\frak a%
}$ and is orthogonal to the subalgebra ${\frak aa}$ as 
\[
ba=bea+l(ba-bea)\theta =l\left( ba\right) \theta ,\quad ab=aeb+l\left(
ab-aeb\right) \theta =l\left( ab\right) \theta 
\]
for all $b\in {\frak b}$, and $baa^{\star }=0=a^{\star }ab$ for all $a\in 
{\frak a}$. And $c$ is an element of a quantum L\'{e}vy B*-algebra ${\frak c}
$, the closure of ${\frak aa}$ in ${\frak a}$ which coincides with all the
combinations $c$ as 
\[
c^{\star }c=a^{\star }a+l\left( c^{\star }c-a^{\star }a\right) \theta ,\quad
\forall a\in {\frak a.} 
\]
Thus $a=b+c$, $bc=0$ for all $a\in {\frak a}$, where $b\in {\frak b}$ is in
a Brownian B*-algebra with the fundamental representation ${\frak b}%
\subseteq {\Bbb C}\oplus {\cal G}_{+}\oplus {\cal G}^{-}$, where ${\cal G}%
_{+}=Pk\left( {\frak a}\right) $, $P=I-E$, ${\cal G}^{-}=k^{*}\left( {\frak a%
}\right) P$, and $c\in {\frak c}$ is in a L\'{e}vy B*-algebra, having the
fundamental representation ${\frak c}\subseteq {\Bbb C\oplus {\cal E}}%
_{+}\oplus {\cal E}^{-}{\Bbb \oplus }{\cal A}$ with non-degenerated operator
algebra ${\cal A}=i\left( {\frak a}\right) $, left and right represented on $%
{\cal E}_{+}=Ek\left( {\frak a}\right) $ and ${\cal E}^{-}=k^{*}\left( 
{\frak a}\right) E$. End of proof.

\end{document}